\documentclass[reqno,12pt]{amsart}
\usepackage{amsfonts}
\usepackage{amssymb}
\usepackage{bbm}
\usepackage{times}
\usepackage{bbm}
\usepackage{amsmath}
\usepackage{txfonts}
\usepackage{stmaryrd}
\usepackage{mathrsfs}

\def\C{\mathbb{C}}

\def\msk{\medskip}

\def\bege{\begin{equation}} \def\ende{\end{equation}}

\def\begr{\begin{eqnarray}} \def\endr{\end{eqnarray}}

\def\bege{\begin{equation}} \def\ende{\end{equation}}
\def\begr{\begin{eqnarray}} \def\endr{\end{eqnarray}}
\def\bnum{\begin{enumerate}} \def\enum{\end{enumerate}}
\def\square{\vbox{
    \hrule height .4pt
    \hbox{\vrule width .4pt height 7pt \kern 7pt
       \vrule width .4pt}
    \hrule height .4pt }}

\begin{document}

\title{On Invariant Subspaces \\ of Subdecomposable Operators}

\author{Junfeng Liu and Songxiao Li}

\address{Junfeng Liu\\ Faculty of Information Technology, Macau University of Science and Technology, Avenida Wai Long, Taipa, Macau.}   \email{jfliu997@sina.com}

\address{Songxiao Li\\ Institute of Fundamental and Frontier Sciences, University of Electronic Science and Technology of China,
610054, Chengdu, Sichuan, P.R. China.  \newline
Institute of Systems Engineering, Macau University of Science and Technology, Avenida Wai Long, Taipa, Macau. } \email{jyulsx@163.com}

\begin{abstract}
In this paper, we prove the Mohebi-Radjabalipour Conjecture under a little additional condition,
and obtain a new invariant subspace theorem for
subdecomposable operators. Our main results contain known results in this topic as special cases.
Moreover, the condition in our theorem is really much weaker,
much simpler and much more nature than that in the known results.
\end{abstract}

\thanks{This project was partially supported by the Macao Science and Technology Development Fund (No.186/2017/A3).}

\maketitle

\section{Introduction and the main result}

The mathematician P. Halmos (\cite{hpr82}, Problem 191) said, ``one of the most important, most difficult, and
most exasperating unsolved problems of operator theory is the problem of
invariant subspaces."

As stated in \cite{avs54}, it was the mathematician and computer scientist Von Neumann who initiated the research work of the invariant
subspace problem in the early thirties of the twentieth century.
More specifically,  Von Neumann proved that every compact operator on a Hilbert space $H$
has a (nontrivial) invariant subspace. In 1954,  N. Aronszajn and  K. Smith \cite{avs54} extended this result from Hilbert spaces to Banach spaces.

We now recall some results which are directly and closely related to the present paper. In 1987,  S. Brown \cite{bro87} showed that if $T$ is a hyponormal operator on a Hilbert space such that the set $\sigma(T)$ is dominating in some nonempty
open set $G$ in the complex plane $\mathbb{C}$, then $T$ has a
(nontrivial) invariant subspace. As pointed out in \cite{mll} and \cite{mor94}, the work initiated by S. Brown \cite{bro87} has been generalized by  J. Eschmeier and  B. Prunaru \cite{esp90}.
More specifically, in 1990,  J. Eschmeier and  B. Prunaru \cite{esp90} proved that if $T$ is a subdecomposable operator on a (general) Banach space such that the set $\sigma(T)$ is dominating in some nonempty open set $G$ in the complex plane $\mathbb{C}$, then $T$ has a
(nontrivial) invariant subspace.

\msk
In 1994, H. Mohebi and M.
Radjabalipour \cite{mor94} raised the following conjecture.

{\bf The Mohebi-Radjabalipour Conjecture} (see \cite{mor94} p.236).
Assume the operators $T\in B(X)$ and $B\in B(Z)$ on Banach spaces
$X$ and $Z$, and the nonempty open set $G$ in the complex plane $\mathbb{C}$
satisfy the following conditions:

(1) \ $qT=Bq$ for some injective $q\in B(X,Z)$ with a closed range $qX$.

(2) \ There exist sequences $\{G(n)\}$ of open sets and $\{M(n)\}$ of invariant
 subspaces of $B$ such that
$\overline{G}(n)\subset G(n+1), G=\cup_nG(n),
\sigma (B|M(n))\subset \mathbb{C}\backslash G(n)$ and
 $\sigma (B/M(n)) \subset \overline{G}(n), n=1,2,\cdots$.

(3) \ $\sigma (T)$ is dominating in $G$.

Then $T$ has a (nontrivial) invariant subspace.

\msk
It is particularly worth mentioning that it can be seen from \cite{mll} (and so on) that
the above conjecture, if true, will contain the main result of \cite{bro78}, \cite{bro87},\cite{esp90} and others as special cases.

\msk

In order to prove the above conjecture, H. Mohebi and M. Radjabalipour \cite{mor94} obtained an important invariant
subspace theorem in reflexive Banach spaces under a additional spectral condition, which is
Theorem I.1 in \cite{mor94} and is also the main result in \cite{mor94}.

Also, in \cite{ml03}, M. Liu proved the above conjecture under the spectral condition that the
set $\sigma(T)\setminus \{\lambda\in \mathbb{C} :~\overline{(\lambda-B^{*})Z^{*}}\neq Z^{*}\} $ is dominating in $G$.
In \cite{mll}, M. Liu and C. Lin proved the above conjecture under another spectral condition that the
set
$$\sigma_0:=\sigma(T)\backslash(\sigma_{p}(B^{**})\cap \{\lambda\in \mathbb{C}: ~\overline{(\lambda-B^{*})\ker q^{*}}
\neq \ker q^{*}\})$$ is dominating in $G$.
Moreover, in \cite{mll}, the authors obtained the result of richness of the invariant subspace Lat$(T)$
for the operator $T$ under the spectral condition that $\sigma_0\cap \sigma_{e}(T)$ is dominating in $G$.

By the way, main results of \cite{ml03} and \cite{mll} contain the main result of \cite{mor94} as a special case (for details, see \cite{ml03} and \cite{mll}).
In particular, the reflexivity condition of the spaces in \cite{mor94} is removed in \cite{ml03} and \cite{mll}.
Moreover, in \cite{mll}, a concrete example is given
to illustrate that there are an operators $A$ that $A$ has
infinitely many invariant subspaces, while its invariant subspace lattice Lat($A$)
is not rich.

\msk

In this paper, we prove the Mohebi-Radjabalipour Conjecture under a little additional condition,
and obtain a new invariant subspace theorem for
subdecomposable operators. In particular,
 we weaken and simplify the spectral conditions in \cite{ml03} and \cite{mll}. More precisely, we obtain
the following invariant subspace theorem.

\msk

\noindent{\bf Theorem 1.}  Assume the operators $T\in B(X)$ and $B\in B(Z)$ on Banach spaces
$X$ and $Z$, and the nonempty open set $G$ in the complex plane $\mathbb{C}$
satisfy the following conditions:

(1) \ $qT=Bq$ for some injective $q\in B(X,Z)$ with a closed range $qX$.

(2) \ There exist sequences $\{G(n)\}$ of open sets and $\{M(n)\}$ of invariant
 subspaces of $B$ such that
$\overline{G}(n)\subset G(n+1), G=\cup_nG(n),
\sigma (B|M(n))\subset \mathbb{C}\backslash G(n)$ and
 $\sigma (B/M(n)) \subset \overline{G}(n), n=1,2,\cdots$.

(3$'$) The set $\sigma (T)\backslash\sigma_r(B^{*})$ is dominating in $G$.\\
Then $T$ has infinitely many invariant subspaces.

Moreover, if the set
$\sigma(T)\backslash(\sigma_p(T^{*})\cup\sigma_p(B^{*})\cup\sigma_r(B^*))$ is
dominating in $G$, then the invariant subspace lattice Lat$(T)$ for the operator
$T$  is rich.

\msk

It can be seen from the definition of the residual spectrum that Theorem 1 contains main results
of \cite{ml03} and \cite{mll} as special cases, and that the condition in Theorem 1 is really much weaker,
much simpler and much more nature than that in \cite{ml03} and \cite{mll}.
Also, since main results of \cite{ml03} and \cite{mll} contain the main result of \cite{mor94} as a special case, Theorem 1 contains main results of
\cite{ml03}, \cite{mll} and \cite{mor94} as special cases.

\section{Preliminary}

In order to prove Theorem 1, we first recall some basic concepts and related results from \cite{esc89}, \cite{esp90}, \cite{esc92}, \cite{ml03} and \cite{mll}.

As usual, we denote by $H^{\infty}(G)$ the Hardy space of all bounded
analytic functions in $G$ equipped with the norm $\|f\|=\mbox{ sup
}\{|f(\lambda)|; \lambda\in G\}$.
We now recall some properties of the Hardy space $H^{\infty}(G)$ from \cite{esc89}.

(a) $H^{\infty}(G)$ is a $w^{*}$-closed subspace of $L^\infty (G)$
with respect to the duality $\langle L^{1}(G), L^{\infty}(G)\rangle$.

(b) A sequence $\{f_n\}$ in $H^\infty (G)$ converges to zero
with respect to the $w^{*}$-topology if and only if it is norm-bounded and
converges to zero uniformly on each compact subset of $G$.

(c) The quotient space $Q:=L^1(G)/^{\bot}H^{\infty}(G)$ is
a separable Banach space.

(d) $H^\infty (G)$ can be identified with the dual space of
the quotient space $Q$.

(e) For every  $\lambda\in G$,
the point evaluation $\mathscr{E}_{\lambda}:H^{\infty}(G)\rightarrow
\mathbb{C},f\rightarrow f(\lambda)$ is a $w^{*}$-continuous linear functional.

For each $f\in H^\infty(G)$ and each $\lambda\in G$, we denote by $f_\lambda$ the unique
function in $H^\infty(G)$ such that
\begin{equation}\label{lamd}
(\lambda-z)f_\lambda(z)=f(\lambda)-f(z)
\end{equation}
for all $z\in G$.

A subset $S$ of the complex plane $\mathbb{C}$ will be called dominating in
the open set $G$ if
$\|f\|=\mbox{sup}\{|f(\lambda)|:\lambda\in G\cap S\}$ holds for
all $f\in H^{\infty} (G)$.

Let $E$ be a nonempty set. Let $n$ and $m$ be positive integers. We define
$$
E^{n}=\{(x_1,x_2,\cdots,x_n):~x_1,x_2,\cdots,x_n\in E\},
$$
$$
M(n\times m,E)=\{(x_{jk})_{1\leq j \leq n,1\leq k\leq m}:~x_{jk}\in E \}.
$$
If $m=n$, we simply write $M(n,E)$ for $M(n\times n, E).$
It is clear that $M(1\times m, E)=E^m.$ Moreover, we
write $M(\infty,E)$ for the set of all infinite matrices
$(x_{jk})_{j,k\geq 1}$ with entries $x_{jk}(j,k=1,2,\cdots)$ in $E$.

If $E$ is a linear space, then $E^n$, $M(n\times m,E)$ and $M(\infty,E)$ are linear spaces
in a natural way. If $E$ is a normed space, then $M(n\times m,E)$ is a normed space and its
norm is defined by
$$\|(x_{jk})_{1\leq j \leq n,1\leq k\leq m}\|=\mbox{max}\{\|x_{jk}\|:j=1,2,\cdots,n, k=1,2,\cdots,m\}$$
for all
$(x_{jk})_{1\leq j \leq n,1\leq k\leq m}\in M(n\times m,E).$
Moreover, if $U$ and $V$ are closed
linear spaces of $E$, then we write

\begin{equation}\label{auv01}
\alpha(U,V)=\inf \{
\|u-v\|;u\in U\mbox{ with }\|u\|=1\mbox{ and }v\in V\}.
\end{equation}

Let $E$ and $F$ be Banach spaces. Then the Banach space of all
bounded linear operators from $E$ into $F$ is denoted by $B(E,F)$. If $F=E$, then we write $B(E)$ for $B(E,E)$.
If $A\in B(E)$, then we denote by $\sigma(A),~\sigma_p(A),~\sigma_r(A)$, $\sigma_e(A)$ and $\rho(A)$
the spectrum, the point spectrum, the residual spectrum, the essential spectrum
and the resolvent set of $A$,  respectively.

The invariant subspace lattice  Lat$(A)$ for an operator $A\in B(E)$ is called to be rich if
there exists an infinite dimensional Banach space $Y$ such that
Lat$(A)$ contains a sublattice order
isomorphic to the lattice Lat$(Y)$ of all
closed linear subspace of $Y$.

\section{Proof of the main result}

Hereafter, we shall assume that $X, Z, T, B, q, G, G(n)$ and $M(n)$ satisfy the conditions (1) and (2) in Theorem 1.
It is easy to show (see also Lemma 1 and  2 in \cite{ml03}) that $q^{*}B^{*}=T^{*}q^{*}$ and $q^{*}\in B(Z^{*},X^{*})$ is surjective.
Moreover, for every positive integer $n$, $M(n)^{\bot}$ is an invariant subspace of $B^{*}$ and
$$
M(n)^{\bot}\subset M(n+1)^{\bot},
\sigma(B^{*}|M(n)^{\bot})\subset \overline{G}(n),
\sigma(B^{*}/M(n)^{\bot})\subset \mathbb{C}\backslash G(n).$$

We write
$
M(G)=\cup_n M(n)^{\bot}.
$
It is easy to see that for any $x\in X,z^{*}\in M(G)$, there is a positive integer
$n$ such that $z^{*}\in M(n)^{\bot}$. Define a
functional $x\otimes z^{*}$ on $H^{\infty}(G)$ as follows:
\begin{equation}\label{xoz01}
x\otimes z^{*}(f)=\langle x,q^{*}f(B^{*}_n)z^{*}\rangle,~f\in H^{\infty}(G),
\end{equation}
where $B_{n}^{*}$ denotes $B^{*}|(M(n)^{\bot})$,
and
$f(B_{n}^{*})$ is defined by the Riesz Functional Calculus.
By the properties of the Riesz Functional Calculus and the $w^*$-topology of the Hardy space
$H^{\infty}(G)$, we can see that
$x\otimes z^{*}$ is a well-defined $w^{*}$-continuous linear
functional on $H^{\infty}(G)$ which is independent of the particular choice of $n$.

Moreover, if $x=(x_1, x_2,\cdots, x_l)\in X^l$ and $z^*=(z^*_{1}, z^*_2,\cdots,z^*_m)\in (M(G))^l$,
then we define
$$x\otimes z^*=(x_j\otimes z^*_k)_{1\leq j \leq l,~~1\leq k\leq m}.$$

\msk
\noindent{\bf Lemma 1}. If $\lambda\in\sigma(T)\backslash(\sigma_p(T^{*})\cup\sigma_p(B^{*})\cup\sigma_r(B^*))$,
then for any finite codimensional subspace $X^{*}_0$ of $X^{*}$, there
exists a sequence $\{z^{*}_n\}$ in $Z^{*}$ such that
$q^{*}z^{*}_n\in X^{*}_{0},\|q^{*}z^{*}_n\|=1$ for every positive integer $n$, and
$$\lim_{n\rightarrow\infty}\|(\lambda-B^{*})z^{*}_n\| =0.$$

{\it Proof}. First of all, we may assume without loss of generality that $\lambda=0$.
Secondly, it is easy to see that $\ker q^{*}$ is an invariant subspace for
$B^{*}$, and so that the quotient operator $B^*/\ker q^*$ induced by $B^*$ on
$Z^*/\ker q^{*}$ is well-defined. Consequently, the operator
$B^*/\ker q^*$ defined by
$$(B^*/\ker q^*)(z^*+\ker q^{*})=B^*z^*+\ker q^{*}, ~~z^*+\ker q^{*}\in Z^*/\ker q^{*}, $$
is a bounded linear operator on $Z^*/\ker q^{*}$.

We now define an operator $\widetilde{q^{*}}$ from
$Z^{*}/\ker q^{*}$ to $X^{*}$ as follows:
$$\widetilde{q^{*}}(z^{*}+\ker q^{*})=q^{*}z^{*}, ~ z^{*}\in Z^{*}.$$
It is not difficult to prove that $\widetilde{q^{*}}$ is a linear bijection.
Thus by the definition of the quotient norm $\|z^{*}+\ker q^{*}\|$
and the Inverse Mapping Theorem, it can be shown that $\widetilde{q^{*}}$
is a (topologically linear) isomorphism from the Banach space $Z^*/\ker q^{*}$ onto the
Banach space $X^*$. Therefore we can assume
without loss of generality that $X^{*}=Z^{*}/\ker q^{*}$.
This implies that $q^{*}$ is the quotient map from $Z^{*}$ onto
$Z^{*}/ \ker q^{*}$, and that $T^{*}$ is the quotient operator induced
by $B^{*}$ on $Z^{*}/ \ker q^{*}$.
That is

\begin{equation}\label{02qz}
  q^*z^*=z^{*}+\ker q^{*},~~z^*\in Z^*,
\end{equation}
and
\begin{equation}\label{02tz}
  T^*(z^{*}+\ker q^{*})=B^*z^*+\ker q^{*}, ~~z^{*}+\ker q^{*}\in Z^{*}/\ker q^{*}.
\end{equation}

Since $X^{*}_0$ is a finite
codimensional subspace of $X^{*}$ , there is a finite dimensional
subspace $X^{*}_{00}$ of $X^{*}$ such that $X^{*}=X^{*}_{0}\oplus X^{*}_{00}$.

We write $Z^{*}_q=\overline{B^{*}(\ker q^{*})} \ (\subset\ker q^{*})$ and define an operator
$\widetilde{B^{*}}$ from $Z^{*}/\ker q^{*}$ to $Z^{*}/Z^{*}_q$ as follows:
\begin{equation}\label{02bz}
\widetilde{B^{*}}(z^{*}+\ker q^{*})=B^{*}z^{*}+Z^{*}_q, ~ z^{*}+\ker q^{*}\in Z^{*}/\ker q^{*}.
\end{equation}
Then it can be verified that $\widetilde{B^{*}}$ is a well-defined bounded
linear operator from
$Z^{*}/\ker q^{*}$ to $Z^{*}/Z^{*}_q$.

We first prove that $B^*$ has dense range, i.e.,
$$\overline{\widetilde{B^{*}}(Z^{*}/\ker q^{*})}=Z^{*}/Z^{*}_q. $$
In fact,
since
$$\lambda=0\in\sigma(T)\backslash (\sigma_p(T^{*})\cup\sigma_p(B^{*})\cup\sigma_r(B^*)),$$
we have $\lambda=0\in \sigma(T), ~\lambda=0\notin \sigma_p(B^*)$
and $\lambda=0\notin \sigma_r(B^*) $. Since
$$\sigma_r(B^*)=\{\lambda\in\mathbb{C}:\lambda-B^*\mbox{ is injective~ and} ~\overline{(\lambda-B^*)Z^*}\neq Z^*\},$$
it follow that
$\overline{B^{*}Z^{*}}=Z^{*}$. This implies that for every vector
$z^{*}+Z^{*}_q\in Z^{*}/Z^{*}_q$ and every real number $\varepsilon>0$, there is a vector $z^{*'} \in
Z^{*}$ such that $\|B^{*}z^{*'}-z^{*}\|<\varepsilon$.
Therefore it follows from (\ref{02bz}) that

\begin{eqnarray*}
  &&\|\widetilde{B^{*}}(z^{*'}+\ker q^{*})-(z^{*}+Z^{*}_q)\|
=\|(B^{*}z^{*'}+Z^{*}_q)-(z^{*}+Z^{*}_q)\|\\
 &=&\|(B^{*}z^{*'}-z^{*})+Z^{*}_q\|
\leq \|B^{*}z^{*'}-z^{*}\|<\varepsilon.
\end{eqnarray*}
This implies $\overline{\widetilde{B^{*}}(Z^{*}/\ker q^{*})}=Z^{*}/Z^{*}_q$.

We now show that ran$\widetilde{B^{*}}$ is not a closed
subspace in the space $Z^{*}/Z^{*}_q$.
Indeed, assume ran$\widetilde{B^{*}}$ is a closed subspace in the space $Z^{*}/Z^{*}_q$. Then
by the above result
we have $\widetilde{B^{*}}(Z^{*}/\ker q^{*})=Z^{*}/Z^{*}_q$. Thus for every $w^{*}+\ker q^{*}\in
Z^{*}/\ker q^{*}$, it follows from the relation $$w^{*}+Z^{*}_q\in
Z^{*}/Z^{*}_q=\widetilde{B^{*}}(Z^{*}/\ker q^{*})$$ that there exists
$z^{*}+\ker q^{*}\in Z^{*}/\ker q^{*}$ such that
$$
w^{*}+Z^{*}_q=\widetilde{B^{*}}(z^{*}+\ker q^{*})=B^{*}z^{*}+Z^{*}_q.
$$
This implies $B^{*}z^{*}-w^{*}\in Z^{*}_q\subset \ker q^{*}.$ Thus by (\ref{02tz}) we have
$$
w^{*}+\ker q^{*}=B^{*}z^{*}+\ker q^{*}=T^{*}(z^{*}+\ker q^{*}).
$$
This implies that $T^{*}$ is a surjection.

 On the other hand, since
$$\lambda=0\in\sigma(T)\backslash (\sigma_p(T^{*})\cup\sigma_p(B^{*})\cup\sigma_r(B^*)),$$
it follows that $\lambda=0\in \sigma(T)$ and $\lambda=0\notin \sigma_p(T^{*})$.
Hence $T^*$ is a injection, and so that $T^{*}$ is a bijection.
Thus by the Inverse Mapping Theorem we have $\lambda=0\in
\rho(T^{*})$, which contradicts $\lambda=0\in \sigma(T)(=\sigma(T^{*}))$.

Since $\widetilde{B^*}$ is a bounded linear operator from $Z^*/\mbox{ker}q^*(=X^*)$
to $Z^*/ Z^{*}_{q}$, it is clear that we can define a bounded linear operator $\overline{B^*}$
from the subspace $X^*_0$ of
$X^{*}$ to $Z^{*}/Z^{*}_q$ as follows:
\begin{equation}\label{02bxw}
  \overline{B^{*}}x^{*}=\widetilde{B^{*}}x^{*}, x^{*}\in X^{*}_0.
\end{equation}

We now show that ran$\overline{B^*}$ is not a closed subspace of the Banach space $Z^{*}/Z^{*}_q$.
In fact, since $X^*=X^*_0\oplus X^*_{00}$, it follows that
$$\widetilde{B^*}X^*=\widetilde{B^*}X^*_0 +\widetilde{B^*} X^*_{00}=\overline{B^*}X^*_0+\widetilde{B^*}X^*_{00}.$$
Since $\widetilde{B^*}X^*_{00}$ is a finite dimensional subspace of $Z^{*}/Z^{*}_q$,
it follows that if $\overline{B^*}X^*_0$ were a closed subspace of $Z^{*}/Z^{*}_q$,
then  $\widetilde{B^*}X^*$ would be a closed subspace of $Z^{*}/Z^{*}_q$, which contradicts that
ran$\widetilde{B^*}$ is not closed.

Then, we show that there is a sequence $\{x^*_n\}$ of unit vectors in $X^*_0$ such that
$\lim_{n\rightarrow \infty}\overline{B^*}x^*_n=0.$

(a). To this end we first show that it is impossible that there is a constant $c>0$ such that
\begin{equation}\label{obx01}
\|\overline{B^*}x^*\|\geq c\|x^*\|
\end{equation}
for all $x^*\in X^*_0$. If not, then it would follows that ran$\overline{B^*}$ is closed
(which contradicts the conclusion in preceding paragraph).
In fact, it suffice to show that for every sequence $\{x^*_n\}$
in $X^*_0$, if $\|\overline{B^*}x^*_n-y^*\|\rightarrow 0$
as $n\rightarrow\infty$, then $y^*\in \mbox{ran}\overline{B^*}$. Indeed,
by (\ref{obx01}), we have
$$\|x^*_n-x^*_m\|\leq \frac{1}{c}\left\|\overline{B^*}(x^*_n-x^*_m)\right\|=\frac{1}{c}\left\|\overline{B^*}x^*_n-\overline{B^*}x^*_m\right\|$$
for all positive integers $n$ and $m$.
Also, since the sequence $\{\overline{B^*}x^*_n\}$ is convergent, the sequence
$\{x^*_n\}$ is a cauchy sequence in $X^*_0$. Since $X^*_0$ is complete, there is
a vector $x^*\in X^*_0$ such that $\|x^*_n-x^*\|\rightarrow 0$
as $n\rightarrow \infty$. Since $\overline{B^*}$ is continuous, it follows that
$\|\overline{B^*}x^*_n-\overline{B^*}x^*\|\rightarrow 0$
as $n\rightarrow \infty$. Also since
$\|\overline{B^*}x^*_n-y^*\|\rightarrow 0$
as $n\rightarrow \infty$, it follows that
$y^*=\overline{B^*}x^*\in \mbox{ran}\overline{B^*}$.

(b). By the conclusion in (a), for each positive integer $n$,
there is an $x^{*'}_n\in X^*_0$ with $x^{*'}_n\neq 0$ such that
$$\|\overline{B^*}x^{*'}_n\|< \frac{1}{n} \left\|x^{*'}_n\right\|, ~~n=1,2,\cdot\cdot\cdot.$$
Set $x^*_n=\frac{1}{\|x^{*'}_n\|}x^{*'}_n, ~~n=1,2,\cdot\cdot\cdot,$ then $\|x^*_n\|=1$
and $\|\overline{B^*}x^{*}_n\|<1/n, ~~n=1,2,\cdot\cdot\cdot.$
This implies that $\lim_{n\rightarrow \infty}\overline{B^*}x^{*}_n=0$.

Since $x^*_n\in X^*_{0}\subset X^*=Z^*/\mbox{ker}q^*$, we can suppose
that $x^{*}_n=u^{*}_n+\ker q^{*}$ for some $u_{_n}^*\in Z^*$. Thus by (\ref{02bz}) and (\ref{02bxw}) we have
$$
\overline{B^*}x_n^*=\widetilde{B^*}x_n^*=B^*u_n^*+Z^{*}_q.
$$
Consequently we can obtain $$\lim_{n\rightarrow\infty}\|B^*u^*_n+Z^{*}_q\|=
\lim_{n\rightarrow\infty}\|\overline{B^{*}}x^{*}_n\|
=0.$$

On the other hand, it follows from the definition of the quotient norm $\|z^*+Z^*_q\|$
that for every positive integer $n$, there exists
$z_n^{'*}\in Z^{*}_q$ such that
$$
\|B^*u^*_n-z_n^{'*}\|<\|B^*u^*_n+Z^{*}_q\|+\frac{1}{n}.
$$
Since $z_{n}^{'*}\in Z^{*}_q=\overline{B^*(\mbox{ker}q^*)}$,
there exists $v_n^*\in \mbox{ker}q^*$ such that
$
\|z_n^{'*}-B^*v_n^*\|<1/n,
$
Consequently we can obtain
$$\lim_{n\rightarrow\infty}\|B^{*}u^{*}_n-B^{*}v^{*}_n\|=0.$$

For every positive integer $n$, we write $z_{_n}^*=u_{_n}^*-v^*_{_n}$,
then we have $z^*_{_n}\in Z^*$, and $\lim_{n\rightarrow\infty}\|B^*z^*_{_n}\|=0$.
Moreover, by (\ref{02qz}) we can obtain
$$q^*z^*_{_n}=q^*u^*_{_n}=u^*_{_n}+\mbox{ker}q^*=x^*_{_n}\in X^*_{_0}.$$
Therefore we have $\|q^*z^*_{_n}\|=\|x^*_{_n}\|=1$ for every positive integer $n$.
The proof is complete.\msk

\noindent{\bf Lemma 2.} Let $m$ be a positive integer. If $\lim_{n\rightarrow
\infty}(\lambda-B^{*})z_{n}^{*}=0$ for some $\lambda \in G(m)$ and
some sequence $\{z_{n}^{*}\}$ in $Z^{*}$, then there is a sequence
$\{\widetilde{{z}_n^{*}}\}$ in $M(m)^{\bot}$ such that
$\lim_{m\rightarrow \infty}(z_n^{*}-\widetilde{{z}_n^{*}})=0$.\msk

{\it \textbf{Proof}.} Since $\lim_{n\rightarrow \infty}(\lambda-B^*)z^*_n=0$, it follows that
\begin{eqnarray*}
   && \|(\lambda-B^*/M(m)^{\bot})(z^*_n+M(m)^{\bot})\|=\|(\lambda-B^*)z^*_n+M(m)^{\bot}\| \\
  &\leq& \|(\lambda-B^*)z^*_n\|\rightarrow 0
\end{eqnarray*}
as $n\rightarrow \infty$, where $\lambda-B^*/M(m)^\perp$ denotes the quotient operator induced by
$\lambda-B^*$ on $Z^*/M(m)^\bot$.
Since $\lambda\in G(m)$ and $\sigma(B^*/M(m)^\bot)\subset \C\setminus G(m)$, it follows that
$\lambda \in \rho(B^*/M(m)^\perp)$, and so that

\begin{eqnarray*}
 &&\inf \{\|z^*_n-z^*\|:z^*\in M(m)^\bot\}=\|z^*_n+M(m)^\bot\|\\
  &=&\|(\lambda-B^*/M(m)^{\bot})^{-1}(\lambda-B^*/M(m)^{\bot})(z^*_n+M(m)^\bot)\|\rightarrow 0
\end{eqnarray*}
as $n\rightarrow \infty.$ This implies that there exists a sequence $\{\widetilde{z^*_n}\}$ in $M(m)^\perp$ such that
$$\lim_{n\rightarrow \infty} \|z^*_n-\widetilde{z^*_n}\|=0.$$
The proof is complete.

\msk

The proof method of Lemmas 3 to 7 below is similar to the corresponding part in
\cite{esc89}, \cite{esp90}, \cite{esc92}, \cite{ml03} and \cite{mll}, but the
context and notations in those paper are different,
thus we still give full proof of these lemmas for the convenience of the reader.
\msk

\noindent{\bf Lemma 3.} Let $r$ and $s$ be positive integers. Let  $c_1,c_2,\cdots,c_r$ be nonnegative
real numbers with $c_1+c_2+\cdots+c_r=1$. If
complex numbers $\lambda_1,\lambda_2,\cdots,\lambda_r \in
(\sigma(T)\backslash(\sigma_p(T^{*})\cup\sigma_p(B^{*})\cup\sigma_r(B^*)))\cap G$, vectors $x_1,x_2,\cdots,x_s\in
X,z_1^{*},z_2^{*},\cdots,$ $z_s^{*}\in M(G)$ and the real number $\varepsilon >0$ are
given, then there are vectors $x\in X$ and $z^{*}\in M(G)$ such that
$\|x\|\leq 3,\|q^{*}z^{*}\|\leq 2$ and

(1) $\|\sum_{k=1}^{r}c_k\mathscr{E}_{\lambda _{k}}-x\otimes z^{*}\|<\varepsilon $;

(2) $\mbox{max}\{\|x\otimes z_{j}^{*} \|:j=1,2,\cdots,s\}
<\varepsilon ,\ \mbox{max}\{\|x_j\otimes z^{*}\|:j=1,2,\cdots,s\}< \varepsilon.$\msk

{\it \textbf{Proof}.} First of all, by the hypothesis, we have
$$\lambda_1,\lambda_2,\cdots,\lambda_r\in
G=\cup_n G(n), ~~z_1^{*},z_2^{*},\cdots,~~~~z_s^{*}\in M(G)=\cup_n
M(n)^{\bot}.$$
Also, since $G(n)\subset G(n+1)$ and $M(n)^{\bot}\subset
M(n+1)^{\bot},n=1,2,\cdots,$ there exists a positive
integer $n$ such that $\lambda_1,\lambda_2,\cdots,\lambda_r\in G(n)
,z_1^{*},z_2^{*},\cdots,z_s^{*}\in M(n)^{\bot}$.

Next, by the $w^{*}$-topological property of $H^{\infty}(G)$,
it can be shown that the set
\begin{equation}
\label{xcq}
X^{*}_c=\{q^{*}f(B_n^{*})z^{*}_j:f\in H^{\infty}(G)\mbox{
with }\|f\|\leq 1,j=1,2,\cdots,s\}
\end{equation}
is a compact subset of $X^*$. Let $\delta$ be an arbitrary
real number with $0<\delta<1$, then there
exist vectors $x^{*}_{-m},\cdots,x^{*}_{-1},x^{*}_0$ in $X^{*}$ such
that the inequality
\begin{equation}
\label{mxx}
\mbox{min}\{\|x^{*}-x^{*}_{-k}\|:k=0,1,\cdots,m\}<\delta
\end{equation}
holds for all $x^{*}\in X^{*}_c$.

By Lemma III.1.1 in \cite{sin81}, there exists a finite codimensional
subspace $X_{1}^{*}$ of $X^{*}$ such that
$${\alpha}({\mbox{span}}
\{x^{*}_{-m},\cdots,x^{*}_{-1},x^{*}_0\},X_{1}^{*})>1-\delta,
$$
where $\alpha$ is defined by (\ref{auv01}),
or equivalently, such that $({\mbox{span}} \{x^*_{-m}, \cdots, x^*_{-1}, x^*_{0}\}) \cap X^{*}_{1}=\{0\}$
and the canonical projection $P$ from the space
$${\mbox{span}} \{x^*_{-m}, \cdots, x^*_{-1}, x^*_{0}\} \oplus X^{*}_{1}$$
onto its subspace ${\mbox{span}} \{x^*_{-m}, \cdots, x^*_{-1}, x^*_{0}\}$ has norm less than
$1/(1-\delta)$. That is,
\begin{equation}\label{pf101}
 \|P\|< \frac{1}{1-\delta}.
\end{equation}
Let
\begin{equation}
\label{x0x}
X^{*}_{0}=\{x_1,x_2,\cdots,x_s\}^{\bot} \subset X^{*}.
\end{equation}
Then $X^*_0$ is a finite codiemensional subspaces of $X^*$, so $X^*_0 \cap X^*_1$ is.
Thus by Lemmas 1 and 2, there are vectors
$x_{1}^{*}\in X_{0}^{*}\cap X_{1}^{*},\widehat{z^{*}_1}\in Z^{*},
\widetilde{z^{*}_1}\in M(n)^{\bot}$
such that
$$
x_{1}^{*}=q^{*}\widehat{z^{*}_1},\|x_{1}^{*}\|=1,\|(\lambda_1-B^{*})\widehat{z^{*}_1}\|
<\delta,\|\widehat{z^{*}_1}-\widetilde{z^{*}_1}\|<\delta.
$$
Again by  Lemma III.1.1 in \cite{sin81}, there is a finite codimensional
subspace $X_{2}^{*}$ in $X^{*}$ such that
$${\alpha}({\mbox{span}}
\{x^{*}_{-m},\cdots,x^{*}_{-1},x^{*}_0,x^{*}_1\},X_{2}^{*})>1-\delta.
$$
Then $X^*_0\cap X^*_1\cap X^*_2$ is a finite codimensional
subspace of $X^*$.
Again by Lemmas 1 and 2, there are vectors
$x_{2}^{*}\in X_{0}^{*}\cap X_{1}^{*}\cap X_{2}^{*},\widehat{z^{*}_2}\in Z^{*},
\widetilde{z^{*}_2}\in M(n)^{\bot}$
such that
$$
x_{2}^{*}=q^{*}\widehat{z^{*}_2},\|x_{2}^{*}\|=1,\|(\lambda_2-B^{*})\widehat{z^{*}_2}\|
<\delta,\|\widehat{z^{*}_2}-\widetilde{z^{*}_2}\|<\delta.
$$
Continuing in this way we can obtain vectors $x^*_k \in X_{0}^{*}\cap X_{1}^{*}\cap\cdots \cap X^*_k$ with $k=1,2,\cdots ,r$, and vectors
$\widehat{z^*_1},\widehat{z^*_2},\cdots,\widehat{z^{*}_r}\in Z^{*},$
$\widetilde{z^*_1},\widetilde{z^*_2},\cdots,\widetilde{z^{*}_r}\in M(n)^{\bot}$
such that
\begin{equation}
\label{xkq}
x_{k}^{*}=q^{*}\widehat{z^{*}_k},\|x_{k}^{*}\|=1,\|(\lambda_k-B^{*})\widehat{z^{*}_k}\|
<\delta,\|\widehat{z^{*}_k}-\widetilde{z^{*}_k}\|<\delta,
\end{equation}
and
\begin{equation}
\label{xmx01}{\alpha}({\mbox{span}}
\{x^{*}_{-m},\cdots,x_{-1}^*, x_{0}^*, x^*_{1},\cdots, x^{*}_{k-1}\},{\mbox{span}}\{x^{*}_k,\cdots,x^{*}_r\})>1-\delta,
\end{equation}
where $k=1,2,\cdots,r.$ Moreover, by (\ref{xmx01}), the inequality
\begin{equation}
\label{xmxk}
{\mbox{max}}\{|\alpha_k|:k=1,2,\cdots,r\}\leq\frac{2}{1-\delta}\left\|\sum^r_{k=1}\alpha_kx_k^*\right\|
\end{equation}
holds for any $\alpha_1,\alpha_2,\cdots,\alpha_r\in \mathbb{ C}$. By (\ref{pf101}) and $x^*_k\in X^*_1$ $(k=1,2,\cdots,r)$,
the canonical projection from the space
$$L={\mbox{span}}\{x^{*}_{-m},\cdots,x^{*}_{-1},x^{*}_{0},x^{*}_{1},\cdots,x^{*}_{r}\}$$
onto its subspace ${\mbox{span}}\{x^{*}_{1},x^{*}_{2},\cdots,x^{*}_{r}\}$
has norm less than $2/(1-\delta)$. Thus by Zenger Lemma
\cite{cz68}, there is a bounded linear functional $l$ on $L$
with $\parallel l\parallel<2/(1-\delta)$ and there are complex numbers
$\mu_1,\mu_2,\cdots,\mu_r$ such that
\begin{equation}\label{kxk}
  \left\|\sum^r_{k=1}\mu_kx_k^*\right\|\leq1, l(\mu_kx_k^*)=c_{_k}\  (k=1,2,\cdots,r),l(x_k^*)=0\  (k=0,-1,\cdots,-m).
\end{equation}
By the canonical isometric isomorphism $L^*\cong X/^{\bot}L$ (notice dim $L<\infty$), there is a
vector $x\in X$ such that $\|x\|<2/(1-\delta)$ and
\begin{equation}\label{lxk}
  l(x_k^*)=\langle x,x_k^*\rangle, k=-m,\cdots,-1,0,1,\cdots,r.
\end{equation}
Hence $\|x\|\leq 3$ if $\delta$ is chosen sufficiently small.

Set $z^*=\sum^r_{k=1}\mu_k\widetilde{z_k^*}$. Then $z^*\in M(G)$. Moreover, it follows from (\ref{xkq}), (\ref{xmxk}) and (\ref{kxk}),
that if $\delta$ is chosen sufficiently small, then
\begin{eqnarray*}
&&\|q^*z^*\|=\left\|\sum^r_{k=1}\mu_kq^*\widetilde{z_k^*}\right\|\leq\left\|\sum^r_{k=1}\mu_kq^*\widehat{z_k^*}\right\|
+\left\|\sum^r_{k=1}\mu_kq^*(\widetilde{z_k^*}-\widehat{z_k^*})\right\|\\
&=&\left\|\sum^r_{k=1}\mu_kx^*_k\right\|+\|q^*\|\delta\sum^r_{k=1}|\mu_k|\leq1+\|q\|\delta r\cdot\left(\frac{2}{1-\delta}\right)\leq 2.
\end{eqnarray*}

Furthermore, it follows from (\ref{lamd}) and (\ref{xkq}) as well as the definition of $B^*_{n}$ and $f(B^*_n)$ that
\begin{eqnarray*}
&&\|(f(\lambda_k)-f(B_n^*))\widetilde{z_k}\|
=\|f_{\lambda_k}(B_n^*)(\lambda_k-B_n^*)\widetilde{z^*_k}\|\\
&\leq&\|f_{\lambda_k}(B^*_n)\|~\|(\lambda_k-B_n^*)\widetilde{z^*_k}\|
=\|f_{\lambda_k}(B^*_n)\|~ \|(\lambda_k-B^*)\widetilde{z^*_k}\|\\
&\leq&\|f_{\lambda_k}(B_n^*)\|(\|(\lambda_k-B^*)\widehat{z_k^*}\|
+\|(\lambda_k-B^*)(\widetilde{z^*_k}-\widehat{z^*_k})\|)\\
&\leq&\|f_{\lambda_k}(B_n^*)\|(\delta+\|\lambda_k-B^*\|\delta),
\end{eqnarray*}
where $f\in H^{\infty}(G)$ and $k=1,2,\cdots,r$. Thus by (\ref{xoz01}), (\ref{xkq}), (\ref{kxk}) and (\ref{lxk}), it can be obtained that
if $\delta$ is chosen sufficiently small, then

\begin{eqnarray}
\label{rxk}
\nonumber &&\left\|\sum^r_{k=1}c_k\mathscr{E}_{\lambda_k}-x\otimes z^*\right\|\\\nonumber
&=&\sup_{f\in H^{\infty}(G), \|f\|\leq 1}~\left|\sum^r_{k=1}c_k\mathscr{E}_{\lambda_k}(f)-x\otimes z^*(f)\right|\\\nonumber
&=&\sup_{f\in H^{\infty}(G), \|f\|\leq 1}~\left|\sum^r_{k=1}\mu_k\langle x,q^*f(\lambda_k)\widehat{z^*_k}\rangle-\sum^r_{k=1}\mu_k\langle x,q^*f(B^*_n)\widetilde{z^*_k}\rangle\right|\\\nonumber
&\leq&\sup_{f\in H^{\infty}(G), \|f\|\leq 1}~\left(\left|\sum^r_{k=1}\mu_k\langle x,q^*f(\lambda_k)(\widehat{z^*_k}-\widetilde{z^*_k})\rangle\right|
+\left|\sum^r_{k=1}\mu_k\langle x,q^*(f(\lambda_k)-f(B^*_n))\widetilde{z^*_k}\rangle\right|\right)\\\nonumber
&\leq&\sup_{f\in H^{\infty}(G), \|f\|\leq 1}~\sum^r_{k=1}|\mu_k|~\|x\|~\|q^*\|(\|f\|\delta
+\|f_{\lambda_k}(B_n^*)\|(\delta+\|\lambda_k-B^*\|\delta))\\
&<&\varepsilon .
\end{eqnarray}

Now it follows from (\ref{xoz01}), (\ref{xcq}),(\ref{mxx}),(\ref{kxk}) and (\ref{lxk})
that if $\delta$ is chosen sufficiently small, then
\begin{eqnarray*}
\|x\otimes z_j^*\|&=&\sup_{f\in H^{\infty}(G), \|f\|\leq 1} |x\otimes z_j^*(f)|=\sup_{f\in H^{\infty}(G), \|f\|\leq 1}|\langle x,q^*f(B^*_n)z^*_j\rangle|\\
&=&\sup_{f\in H^{\infty}(G), \|f\|\leq 1}~~\min_{0\leq k\leq m} |\langle x,q^*f(B^*_n)z^*_{j}-x^*_{-k}\rangle|\\
&\leq&\sup_{f\in H^{\infty}(G), \|f\|\leq 1}~~\min_{0\leq k\leq m} \|x\|~\| q^*f(B^*_n)z^*_{j}-x^*_{-k}\|\leq\|x\|\delta\\
&<&\varepsilon,
\end{eqnarray*}
where $j=1,2,\cdots,s$.

Finally, since $x^*_k\in X^*_0~~(k=1,2,\cdots,r)$, it follows from
(\ref{x0x}) and (\ref{xkq}) that
$$\langle x_j,q^*\widehat{z^*_k}\rangle=\langle x_j,x^*_k\rangle=0,$$
where $j=1,2,\cdots, s$ and $k=1,2,\cdots,r$. Thus,
if $\delta$ is chosen sufficiently small, then by the same calculation as in (\ref{rxk}), we get
\begin{eqnarray*}
\|x_j\otimes z^*\|&=&\sup_{f\in H^{\infty}(G), \|f\|\leq 1}|x_j\otimes z^*(f)|
=\sup_{f\in H^{\infty}(G), \|f\|\leq 1}\left|\sum^r_{k=1}\mu_k\langle x_j,q^*f(B_n^*)\widetilde{z^*_k}\rangle \right|\\
&=&\sup_{f\in H^{\infty}(G), \|f\|\leq 1} \left|\sum^r_{k=1}\mu_k\langle x_j,q^*f(\lambda_k)\widehat{z^*_k}\rangle-\sum^r_{k=1}\mu_k\langle x_j,q^*f(B_n^*)\widetilde{z^*_k}\rangle\right|\\
&<& \varepsilon,
\end{eqnarray*}
where $j=1,2,\cdots, s$.
The proof is complete.

\msk

\noindent{\bf Lemma 4.} Let $r$ and $s$ be positive integers. Let  $c_1,c_2,\cdots,c_r$ be complex numbers with $|c_1|+|c_2|+\cdots+|c_r|\leq 1$. If
complex numbers $\lambda_1,\lambda_2,\cdots,\lambda_r \in
(\sigma(T)\backslash(\sigma_p(T^{*})\cup\sigma_p(B^{*})\cup\sigma_r(B^*)))\cap G$, vectors $x_1,x_2,\cdots,x_s\in
X,z_1^{*},z_2^{*},\cdots,$ $z_s^{*}\in M(G)$ and the real number $\varepsilon >0$ are
given, then there are vectors $x\in X$ and $z^{*}\in M(G)$ such that $\|x\|\leq 12,  \|q^*z^*\|\leq 8$ and

(1) \ \ $\|\sum_{k=1}^{r}c_k\mathscr{E}_{\lambda _{k}}-x\otimes z^{*}\|<\varepsilon $;

(2) \ \ $\mbox{ max }\{\|x\otimes z_{j}^{*} \|:j=1,2,\cdots,s\}
<\varepsilon ,\ \mbox{max}\{\|x_j\otimes z^{*}\|:j=1,2,\cdots,s\}< \varepsilon.$

{\it \textbf{Proof}.} \ \ Write $c_k=c_{1k}-c_{2k}+ic_{3k}-ic_{4k}$ with $c_{mk}\geq 0$,
$k=1,2,\cdots,r,~~m=1,2,3,4,~~i=\sqrt{-1}$. Then $0\leq c_{mk}\leq |c_k|$
and
$$
\sum_{k=1}^{r}c_{mk}\leq\sum_{k=1}^{r}|c_k|\leq 1,\ \ m=1,2,3,4.
$$

For each given $m=1,2,3,4$, if $\sum^r_{k=1}c_{mk}=0$, then we take $c^{'}_{mk}=0\ (k=1,2,\cdots,r), \ x^{(m)}=0, \ z^{*(m)}=0.$
If $\sum^r_{k=1}c_{mk}\neq 0$, then we put $c_{mk}^{'}=c_{mk}/\sum_{k=1}^{r}c_{mk}$ for all $k=1,2,\cdots,r.$
Consequently $c_{mk}^{'}\geq 0,$ and $c_{m1}^{'}+c_{m2}^{'}+\cdots+c_{mr}^{'}=1.$
By Lemma 3, we can choose successively pairs
$(x^{(m)},z^{*(m)})_{1\leq m\leq 4}$ of vectors
$x^{(m)}\in X$ and $z^{*(m)}\in M(G)$ such that
$\|x^{(m)}\|\leq 3,\|q^*z^{*(m)}\|\leq2$, and
$$
\left\|\sum_{k=1}^r c^{'}_{mk}\mathscr{E}_{\lambda_k}-x^{(m)}\otimes z^{*(m)}\right\|
<\frac{\varepsilon}{16},$$
$$
\|x^{(m)}\otimes z^*_j\|<\frac{\varepsilon}{16},~
\|x_j\otimes z^{*(m)}\|<\frac{\varepsilon}{16},\
\|x^{(m)}\otimes z^{*(l)}\|<\frac{\varepsilon}{16}
$$
for every $j=1,2,\cdots,s,$ and all $m,l=1,2,3,4$ with $l\not=m.$

Set
\begin{eqnarray*}
x&=&\left(\sum_{k=1}^{r}c_{1k}\right)^{1/2}\ x^{(1)}
-\left(\sum_{k=1}^{r}c_{2k}\right)^{1/2}\ x^{(2)}
+i\left(\sum_{k=1}^{r}c_{3k}\right)^{1/2}\ x^{(3)}
-i\left(\sum_{k=1}^{r}c_{4k}\right)^{1/2}\ x^{(4)},\\
z^*&=&\left(\sum_{k=1}^{r}c_{1k}\right)^{1/2}\ z^{*(1)}
+\left(\sum_{k=1}^{r}c_{2k}\right)^{1/2}\ z^{*(2)}
+\left(\sum_{k=1}^{r}c_{3k}\right)^{1/2}\ z^{*(3)}
+\left(\sum_{k=1}^{r}c_{4k}\right)^{1/2}\ z^{*(4)}.
\end{eqnarray*}
Then $\|x\|\leq 12,\ \|qz^*\|\leq 8$ and
\begin{eqnarray*}
&&\left\|\sum^r_{k=1} c_k\mathscr{E}_{\lambda_k}-x\otimes z^*\right\|\\
&\leq&\left(\sum_{k=1}^{r}|c_k|\right)\left(\sum_{m=1}^{4}\left\|\sum_{k=1}^{r}c_{mk}^{'}\mathscr{E}_{\lambda_{k}}-x^{(m)}\otimes z^{*(m)}\right\|+
\sum_{1\leq m,l\leq 4, l\neq m}\|x^{(m)}\otimes z^{*(l)}\|\right)<\varepsilon,
\end{eqnarray*}
and the inequalities
$$
\|x\otimes z^*_j\|\leq\left(\sum_{k=1}^{r}|c_k|\right)^{1/2}
\left(\sum_{m=1}^{4}\|x^{(m)}\otimes z^*_{j}\|\right)<\varepsilon,
$$
and
$$
\|x_j\otimes z^*\|\leq\left(\sum_{k=1}^{r}|c_k|\right)^{1/2}
\left(\sum_{m=1}^{4}\|x_j\otimes z^{*(m)}\|\right)<\varepsilon
$$
hold for every $j=1,2,\cdots,s$.
The proof is complete.

\msk

We now denote by $Q_0$ the set of those vectors $L$ in $Q=L^{1}(G)/^{\perp} H^{\infty}(G)$
with the following property:

For any given positive integer $s$,
vectors
$$x_1,x_2,\cdots,x_s\in X,
z^*_1, z^*_2, \cdots,\ \ \ \ z^*_s\in M(G)$$
and the real number
$\varepsilon>0$, there are vectors $x\in X$ and $z^{*}\in M(G)$ such that
$\|x\|\leq 1,\|q^{*}z^{*}\|\leq 1,
\|L-x\otimes z^{*}\|<\varepsilon$ and
$$
\mbox{ max }\{\|x\otimes z_{j}^{*} \|:j=1,2,\cdots,s\}
<\varepsilon ,\ \mbox{max}\{\|x_j\otimes z^{*}\|:j=1,2,\cdots,s\}< \varepsilon.
$$

It is easy to show that the space $Q_{0}$ is a (norm) closed linear subspace of $Q$.

\msk

\noindent{\bf Lemma 5.} If the set
$\sigma(T)\backslash(\sigma_p(T^{*})\cup\sigma_p(B^{*})\cup\sigma_r(B^*))$ is
dominating in $G$, then $\{L\in Q:\|L\|\leq 1/96\}\subset Q_0$.

{\it \textbf{Proof}.} By Lemma 4 the absolutely-convex of the set
$$\left\{ \frac{1}{96}\mathscr{E}_\lambda\ :\lambda\in(\sigma(T)\backslash(\sigma_p(T^{*})\cup\sigma_p(B^{*})\cup\sigma_r(B^*)))\cap G\right\}$$
is contained in $Q_0$, that is
\begin{equation}\label{ave}
 \mbox{aco} \left\{ \frac{1}{96}\mathscr{E}_\lambda\ :\lambda\in(\sigma(T)\backslash(\sigma_p(T^{*})\cup\sigma_p(B^{*})\cup\sigma_r(B^*)))\cap G\right\}\subset Q_0.
\end{equation}
Since the set $\sigma(T)\backslash(\sigma_p(T^{*})\cup\sigma_p(B^{*})\cup\sigma_r(B^*))$ is dominating in $G$,
it follows from Proposition 2.8 in \cite{bcp79} that
\begin{equation}\label{mac}
  \overline{\mbox{aco}} \{\mathscr{E}_\lambda\ :\lambda\in(\sigma(T)\backslash(\sigma_p(T^{*})\cup\sigma_p(B^{*})\cup\sigma_r(B^*)))\cap G\}=\{L\in Q :\|L\|\leq1\}.
\end{equation}
Since $Q_0$ is a closed subset of $Q$, it follows from (\ref{ave}) and (\ref{mac}) that
\begin{eqnarray*}
   &&\left\{L\in Q :\|L\|\leq \frac{1}{96}\right\} \\
  &=&\overline{\mbox{aco}}\left\{ \frac{1}{96}\mathscr{E}_\lambda :
\lambda\in(\sigma(T)\backslash(\sigma_p(T^{*})\cup\sigma_p(B^{*})\cup\sigma_r(B^*)))\cap G\right\}\subset Q_0.
\end{eqnarray*}
The proof is complete.

\msk

\noindent{\bf Lemma 6.} Let $n$ be a positive integer. Let vectors $M=(M_{jk})_{1\leq j,k\leq n}\in M(n,Q)$,
$x_0 \in X^n, z_0^*\in (M(G))^n$ and the real number $\varepsilon>0$ be given.
If the set $\sigma(T)\backslash(\sigma_p(T^{*})\cup\sigma_p(B^{*})\cup\sigma_r(B^*))$ is dominating in $G$,
then there are vectors $x\in X^n$ and $z^*\in (M(G))^n$ such that

(1)$\|M-x\otimes z^*\|<\varepsilon$;

(2)$\max(\|x-x_0\|,\|q^*(z^*-z_0^*)\|)\leq10n\|M-x_0\otimes z_0^*)\|^{1/2}$,\\
where $x\otimes z^*=(x_j\otimes z^*_k)_{1\leq j, k\leq n}$ and $q^* z^*=(q^* z^*_1, q^* z^*_2,\cdots, q^*z^*_n)$
whenever
$x=(x_1, x_2,\cdots, x_n)\in X^n$ and $z^*=(z^*_1, z^*_2, \cdots, z^*_n)\in (M(G))^n$.

{\it \textbf{Proof}.} Set $c=\parallel M-x_0\otimes z_0^*\parallel$.
If $c=0$, then we take $x=x_0$, $z^*=z^*_0$.
If $c>0$, then we write
$$L=\left(\frac{1}{100n^2c}(M-x_0\otimes z^*_0)\right)=(L_{jk})_{1\leq j,k\leq n}.$$
Therefore $\|L\|=1/(100n^2)$, this implies
$\|n^2L_{jk}\|\leq1/100$ for all $j,k=1,2,\cdots,n$. Thus by Lemma 5, we have
$n^2L_{jk}\in Q_0$ for all $j,k=1,2,\cdots,n$.
Let $\delta$ be an arbitrary real number with $0<\delta<1$.
Then by the definition of $Q_0$, we
can choose successively (in an arbitrary order) pairs
$(x_{jk},z^*_{jk})_{1\leq j,k\leq n}$
of vectors $x_{jk}\in X$ and $z^*_{jk}\in M(G)$
such that max$(\|x_{jk}\|, \|q^*z^*_{jk}\|)\leq 1$ and
$$ \|n^2L_{jk}-x_{jk}\otimes z^*_{jk}\|<\delta,\ \ \ \ \ \ \|x_{jk}\otimes z^*_{lm}\|<\delta, $$
\begin{equation}\label{xjjk}
\|x_{jk}\otimes z^*_0\|<\delta,\ \ \ \ \ \ \|x_0\otimes z^*_{jk}\|<\delta,
\end{equation}
where $j,k,l,m=1,2,\cdots,n$ with $(l,m)\neq (j,k).$
Write
$$x_j=\frac{1}{n}\sum^n_{l=1}x_{jl},\ \ z_k^*=\frac{1}{n}\sum^n_{m=1}z^*_{mk},j,k=1,2,\cdots,n.$$
$$\widetilde{x}=(x_1,x_2,\cdots,x_n),\ \ \widetilde{z^*}=(z^*_1,z^*_2,\cdots,z^*_n).$$
Then it follows from (\ref{xjjk}) that
\begin{eqnarray*}
&&\|L_{jk}-x_j\otimes z^*_{k}\|\\
&=&\left\|\frac{1}{n^2}\left(n^2L_{jk}-\sum^n_{l,m=1}x_{jl}\otimes z^*_{mk}\right)\right\| \\
&\leq&\frac{1}{n^2}\left(\|n^2L_{jk}-x_{jk}\otimes z_{jk}^*\|+\sum_{(j,l)\neq(m,k)}\left\|x_{jl}\otimes z^*_{mk}\right\|\right) \\
&<&\delta,
\end{eqnarray*}
where $j,k=1,2,\cdots,n$, so that
\begin{equation}\label{121}
  \left\|\frac{1}{100n^2c}(M-x_0\otimes z_0^*)-\widetilde{x}\otimes\widetilde{z^*}\right\|
  =\|L-\widetilde{x}\otimes\widetilde{z^*}\|=\max_{1\leq j,k\leq n}\|L_{jk}-x_j\otimes z^*_{k}\|<\delta.
\end{equation}
Moreover, it follows from (\ref{xjjk}) that
\begin{equation}\label{z8x}
  \|\widetilde{x}\otimes z^*_0\|<\delta, \ \ \ \|x_0\otimes \widetilde{z^*}\|<\delta,
\end{equation}
and
$\|\widetilde{x}\|\leq 1,~~\|q^*\widetilde{z^*}\|\leq 1.$

Define $x=x_0+10nc^{^{1/2}}\widetilde{x}$,
$z^*=z^*_0+10nc^{^{1/2}}\widetilde{z^*}$.
Then we have
$$\|x-x_0\|=10nc^{^{1/2}}\|\widetilde{x}\|\leq10n\|M-x_0\otimes z^*_0\|^{^{1/2}},$$
and
$$\|q^*(z^*-z^*_0)\|=10nc^{^{1/2}}\|q^*\widetilde{z^*}\|\leq10n\|M-x_0\otimes z^*_0\|^{^{1/2}}.$$
Moreover by (\ref{121}) and (\ref{z8x}) we can see that the inequality
\begin{eqnarray*}
&&\|M-x\otimes z^*\|\\
&\leq& \|M-x_0\otimes z^*_0-100n^2c\widetilde{x}\otimes \widetilde{z^*}\|
+10nc^{^{1/2}}(\|x_0\otimes \widetilde{z^*}\|+\|\widetilde{x}\otimes z^*_0\|)\\
&\leq& 100n^2c\left\|\frac{1}{100n^2c}(M-x_0\otimes z^*_0)-\widetilde{x}\otimes \widetilde{z^*}\right\|+20nc^{1/2}\delta\\
&<& \varepsilon
\end{eqnarray*}
holds if $\delta$ is chosen sufficiently small.
The proof is complete.

\msk

\noindent{\bf Lemma 7}. If the set $\sigma(T)\backslash(\sigma_p(T^{*})\cup\sigma_p(B^{*})\cup\sigma_r(B^*))$ is
dominating in $G$, then for each infinite matrix $L=(L_{jk})_{j,k\geq 1}\in
M(\infty,Q)$, there are sequences $\{x_n\}$ and $\{z^*_n\}$ such
that

(1) \ \ $x_n\in X^{n}$ and $z^*_n \in (M(G))^{n}$;

(2) \ \ for each positive integer $j$, the limits
$x(j)=\lim_{n\rightarrow \infty}x_{n}(j)\in X$ and
$x^{*}(j)=\lim_{n\rightarrow \infty}q^{*}z^*_n(j) \in
X^{*}
$
exist, where $x_n(j)$ and $z_{n}^{*}(j)$ denote the $j$th components
of $x_n$ and $z_{n}^{*}$, respectively;

(3) \ \  for all positive integers $j$ and $k$, we have
$
L_{jk}=\lim_{n\rightarrow\infty}x_n(j)\otimes z_{n}^{*}(k),
$
where the limit is taken in $Q$, while $Q=L^{1}(G)/^{\bot}H^{\infty}(G).$
\msk

{\it \textbf{Proof}.} For every positive integer $n$, write
$$\alpha_n=\frac{1}{n^{6} (\mbox{max}\{\|L_{jk}\|: j,k=1,2,\cdots, n\}+1)},$$
$$M=(\alpha_j\alpha_kL_{jk})_{j,k\geq 1}=(M_{jk})_{j,k\geq 1}\in M(\infty, Q).$$
Take $x_0=0\in X^1,\ z^*_0=0\in (M(G))^1$.
By Lemma 6 and the definition of $\alpha_1$,
there are vectors $x_1\in X^1$ and $\ z^*_1\in (M(G))^1$ such that
\begin{eqnarray*}
(\mbox{a}) ~~~~~~\|M_{11}-x_1\otimes z^*_1\|<1/2^6,\ \ \ \ \ \ \ \ \ \ \ \ \ \ \ \ \ \ \ \ \ \ \ \ \ \ \ \ \ \ \ \ \ \ \ \ \ \ \ \ \ \
\end{eqnarray*}
\begin{eqnarray*}
(\mbox{b})~&&\mbox{max}(\|x_1-x_0\|, \ \|q^*(z_1^*-z^*_0)\|)\leq10\|M_{11}-0\|\ \ \ \ \ \ \ \\
&&=10\alpha_1\alpha_1\|L_{11}\|\leq10\alpha_1\|L_{11}\|<10.
\end{eqnarray*}

By induction, assume that there are vectors $x_{n-1}\in X^{n-1}$ and $\ z^*_{n-1}\in M(G)^{n-1}$ such that
\begin{equation}\label{mjk}
(\mbox{a}')~~ \|(M_{jk})_{1\leq j,k\leq n-1}-x_{n-1}\otimes z^*_{n-1}\|<\frac{1}{n^6},\ \ \ \ \ \ \ \ \ \ \ \ \ \ \ \ \ \ \ \ \ \ \ \ \ \ \ \ \
\end{equation}
\begin{eqnarray*}
(\mbox{b}')~~\mbox{max}(\|x_{n-1}-x_{n-2}\|,\|q^*(z^*_{n-1}-z^*_{n-2})\|)\leq \frac{10}{(n-1)^2}.~~~~~~~~~~~~
\end{eqnarray*}
Then by Lemma 6, there are vectors $x_n\in X^n$ and $\ z^*_n\in M(G)^n$
such that
\begin{equation}\label{mjk2}
(\mbox{a}'')\ \|(M_{jk})_{1\leq j,k\leq n}-x_{n}\otimes z^*_{n}\|<\frac{1}{(n+1)^6},\ \ \ \ \ \ \ \ \ \ \ \ \ \ \ \ \ \ \ \ \ \ \ \ \ \ \ \ \
\end{equation}
\begin{equation}\label{mjk3}
(\mbox{b}'') \mbox{max}(\|x_{n}-x_{n-1}\|,\|q^*(z^*_{n}-z^*_{n-1})\|)\leq 10n\|(M_{jk})_{1\leq j,k\leq n}-x_{n-1}\otimes z^*_{n-1}\|^{1/2},
\end{equation}
where
$$x_{n-1}=(x_{n-1}(1),x_{n-1}(2),\cdots,x_{n-1}(n-1))\cong(x_{n-1}(1),x_{n-1}(2),\cdots,x_{n-1}(n-1),0)\in X^n,$$
$$z^*_{n-1}=(z^*_{n-1}(1), z^*_{n-1}(2),\cdots,z^*_{n-1}(n-1)) \cong (z^*_{n-1}(1), z^*_{n-1}(2),\cdots,z^*_{n-1}(n-1),0)\in (M(G))^n.$$

Moreover, by the definition of $\alpha_n$ we can obtain
$$\|M_{jn}-0\|=\alpha_j\alpha_n\|L_{jn}\|< \alpha_n\|L_{jn}\|<\frac{1}{n^6},$$
$$\|M_{nj}-0\|=\alpha_n\alpha_j\|L_{nj}\|< \alpha_n\|L_{nj}\|<\frac{1}{n^6},$$
where $j=1,2,\cdots,n.$ Thus by (\ref{mjk}) and (\ref{mjk3}) we have
\begin{eqnarray*}
  &&\|x_{n}-x_{n-1}\| \leq 10n \|(M_{jk})_{1\leq j,k\leq n}-x_{n-1}\otimes z^*_{n-1}\|^{1/2} \\
  &=& 10n(\max(\|(M_{jn})_{1\leq j,k\leq n-1}-x_{n-1}\otimes z^*_{n-1}\|), \max_{1\leq j\leq n}(\|M_{jn}-0\|, \|M_{nj}-0\|))^{1/2}  \\
  &<&  \frac{10}{n^2}.
\end{eqnarray*}
Similarly, $\|q^*(z^*_n-z^*_{n-1})\|<10/n^2.$
Thus for each positive integer $j$, we have
$$\|x_{n}(j)-x_{n-1}(j)\|\leq \|x_n-x_{n-1}\|<\frac{10}{n^2}$$
and
$$\|q^*z^*_n(j)-q^*z^*_{n-1}(j)\|\leq \|q^*(z^*_n-z^*_{n-1})\|<\frac{10}{n^2},$$
where $n=1,2,\cdots$.
This implies that $\{x_n(j)\}$ and $\{q^*z^*_n(j)\}$ are Cauchy sequences in Banach
spaces $X$ and $X^*$, respectively. Also, since $X$ and $X^*$ are Banach spaces,
there are vectors $x(j)\in X$ and $x^*(j)\in X^*$
such that
\begin{equation}\label{lnr01}
\lim_{n\rightarrow \infty}x_n(j)=x(j),\ \lim_{n\rightarrow\infty}q^*z_n^*(j)=x^*(j),
\end{equation}
where $x_n(j)$, $x(j)$, $z_n^*(j)$ and $x^*(j)$ denote the $j$th
components of $x_n$, $x$, $z^*_n$ and $z^*$, respectively.

On the other hand, it follows from
(\ref{mjk2}) that
$$\|(\alpha_{j}\alpha_{k}L_{jk})_{1\leq j,k\leq n}-(x_{n}(j)\otimes z^*_{n}(k))_{1\leq j,k\leq n}\|<\frac{1}{(n+1)^6}.$$
So
 $$ \|L_{jk}-(x_n(j)/\alpha_j)\otimes (z_n^*(k)/\alpha_k)\|< \frac{1}{\alpha_j\alpha_k}\cdot\frac{1}{(n+1)^6}\rightarrow 0$$
as $n\rightarrow \infty$,
where $j,k=1,2,\cdots, n.$ That is,

\begin{equation}\label{lnr003}
  \lim_{n\rightarrow \infty}(x_n(j)/\alpha_j)\otimes (z^*_n(k)/\alpha_k)=L_{jk}
\end{equation}
for all positive integers $j$ and $k$, where the limit is taken in $Q$.

In order to obtain the conclusion,
we write $$\widetilde{x_n}(j)=x_n(j)/\alpha_{j},~\widetilde{z^*_n}(j)=z^*_n(j)/\alpha_{j},
~\widetilde{x}(j)=x(j)/\alpha_{j},~\widetilde{z^*}(j)=z^*(j)/\alpha_{j}
$$
and
$$\widetilde{x_n}=(\widetilde{x_n}(1), \widetilde{x_n}(2), \cdots, \widetilde{x_n}(n)),~~
\widetilde{z^*_n}=\widetilde{z^*}(1), ~~\widetilde{z^*}(2),\cdots,\widetilde{z^*}(n).$$
Then by (\ref{lnr01}) and (\ref{lnr003}), we have

(1) $\widetilde{x_n}\in X^n$ and $\widetilde{z^*_n}\in (M(G))^n, n=1,2,\cdots,$

(2) $\widetilde{x}(j)=x(j)/\alpha_{j}=\lim_{n\rightarrow \infty}x_{n}(j)/\alpha_j=\lim _{n\rightarrow \infty }
\widetilde{x_{n}}(j),~j=1,2,\cdots,
$
and
$\widetilde{x^*}(j)=x^*(j)/\alpha_j=\lim_{n\rightarrow \infty}(q^*z^*_{n}(j))/\alpha_j=\lim_{n\rightarrow \infty}q^*\widetilde{z^*}_n(j), ~~j=1,2,\cdots,$

(3) $L_{jk}=\lim_{n\rightarrow \infty}(x_{n}(j)/\alpha_j)\otimes (z^*_n(k)/\alpha_k)=\lim_{n\rightarrow \infty}
\widetilde{x_n}(j)\otimes \widetilde{z^*_n}(k), ~j,k=1,2,\cdots.$

The proof is complete.\msk

We are now in a position to prove the main result of this paper.\msk

\textbf{ {\it Proof of Theorem 1.}} The proof is divided into two steps.

{\bf Step 1}.  we first prove that if $X,~Z,~T,~B,~q,~G,~G(n),$ and $M(n)$ satisfy
the conditions (1) and (2)
in Theorem 1, and if the set $\sigma(T)\backslash(\sigma_p(T^{*})\cup\sigma_p(B^{*})\cup\sigma_r(B^*))$
is dominating in $G$, then
the invariant subspace lattice Lat$(T)$ for the operator
 $T$  is rich.
By means of Lemma 7, the proof can be completed. Moreover, the specific structure
of the invariant subspace lattice Lat$(T)$ can be seen from the proof.
To this end, fix a
complex number $\lambda \in G$. Since the set
$\sigma(T)\backslash(\sigma_p(T^{*})\cup\sigma_p(B^{*})\cup\sigma_r(B^*))$ is
dominating in $G$, it follows from Lemma 7 that there are sequences
$\{x_n\},\{z_{n}^{*}\},\{x(j)\}$ and $\{x^{*}(k)\}$ which satisfy (1),(2) and (3) in Lemma 7 with respect to the matrix
$$L=(L_{jk})_{j,k\geq 1}=(\delta_{jk}\mathscr{E}_{\lambda})_{j,k \geq 1}\in
M(\infty,Q),$$  where $\delta_{jk}$ denotes the Kronecker delta. Since
$z_{n}^*(k)\in M(G)$ for all
positive integers $n$ and $k$, it follows that
for any given positive integers $n$ and $k$, there is a positive integer
 $m=m(n,k)$ such that
 $z_{n}^{*}(k)\in
M(m(n,k))^{\bot}$. Consequently for each polynomial $p \in C[z]$ in
one complex variable, and all positive integers $j,k$,  by the relation
$q^*B^*=T^*q^*$, we get
\begin{eqnarray}\label{02jkp}
\delta_{jk}p(\lambda)&=&\delta_{jk}\mathscr{E}_{\lambda}(p)=\lim
_{n\rightarrow\infty}\langle
x_n(j),q^*p(B^*|M(m(n,k))^{\bot})z^*_n(k)\rangle\nonumber\\
&=&\lim_{n\rightarrow\infty}\langle x_n(j),
p(T^{*})q^*z^*_n(k)\rangle
=\langle x(j), p(T^{*})x^{*}(k)\rangle \nonumber\\
&=&\langle p(T)x(j),
x^{*}(k)\rangle.
\end{eqnarray}
We now define two subspaces in $X$ as follows:
$$
U=\overline{\mbox{span}}\{p(T)x(j):~p\in C[z],j=1,2,\cdots\},
$$
$$
V=\overline{\mbox{span}}\{p(T)(\lambda-T)x(j):~p\in
C[z],j=1,2,\cdots\}.
$$
It is easy to see that $U,V\in \mbox{Lat}(T)$, and $V\subset U$.
Moreover, it follows from (\ref{02jkp}) that the equality
$$\langle y,x^*(k)\rangle=\langle p(T)(\lambda-T)x(j),x^*(k)\rangle=\delta_{_{jk}}p(\lambda)(\lambda-\lambda)=0$$
holds for every vector of the form $y=p(T)(\lambda-T)x(j)$ in $V$ and every $k=1,2,\cdot\cdot\cdot.$
This implies that the equality $\langle y,x^*(k)\rangle=0$ holds for every $y\in V$ and every $k=1,2,\cdot\cdot\cdot.$
Hence $$(U/V)^*=V^{\perp}\supset \mbox{span}\{x^*(k)|k=1,2,\cdot\cdot\cdot\}.$$
On the other hand, by (\ref{02jkp}), $\langle x(j), x^*(k)\rangle=\delta_{_{jk}}$.
This implies that $\{x^*(k)|k=1,2,\cdot\cdot\cdot\}$
is a linearly independent subset in $X^*$. Consequently,
$$\mbox{dim}(U/V)^*\geq \mbox{dim}(\mbox{span}\{x^*(k)|k=1,2,\cdot\cdot\cdot\})=\infty.$$
Therefore
$U/V$ is a infinite dimensional Banach space. Let
$\pi:U\rightarrow U/V$ be the quotient map. Then
the map $$\psi:\mbox{Lat}(U/V)\rightarrow \mbox{Lat}(T),
W\rightarrow \pi^{-1}(W)$$ is a lattice embedding, where Lat$(U/V)$ denotes the
the lattice of all closed linear subspaces of the Banach space $U/V$.
Consequently the invariant subspace lattice Lat$(T)$ for the operator $T$  is rich.

{\bf Step 2.} Now we show that if  $X,~Z,~T,~B,~q,~G,~G(n),$ and $M(n)$ satisfy
the conditions (1), (2) and (3$'$) in Theorem 1, then $T$ has infinitely many invariant subspaces.
Indeed, by the spectral theory, we have
$$\sigma_p(T^*)\subset\sigma_p(T)\cup \sigma_r(T).$$
Therefore if  there is a complex number $\lambda\in \sigma_p(T^*)$, then we have $\lambda\in \sigma_p(T)\cup \sigma_r(T)$.
This implies that $\lambda\in \sigma_p(T)$ or $\lambda\in \sigma_r(T)$.
If $\lambda\in \sigma_p(T)$, then $M_{_\lambda}:=$ker$(\lambda-T)$
is a (nontrivial) invariant subspace for $T$.   If $\lambda\in \sigma_r(T)$,
then $M_{_\lambda}:=\overline{ \mbox{ran}(\lambda-T)}(\neq X)$ is a
(nontrivial) invariant subspace for $T$.
Consequently if $\sigma_p(T^*)$ is a infinite set,
then $T$ has a infinitely many invariant subspaces.
By using the condition (1) in Theorem 1, it is easy to show in a similar
way that if $\sigma_p(B^*)$ is a infinite set,
then $T$ has a infinitely invariant subspaces.
According to the above result, we can assume that  $\sigma_p(T^*)\cup\sigma_p(B^*)$ is a finite set. Since
$\sigma(T)\backslash\sigma_r(B^{*})$  is dominating in $G$, it
follows from the maximum modulus principle for the analytic function
that the set $\sigma(T)\backslash(\sigma_p(T^{*})\cup\sigma_p(B^{*})\cup\sigma_r(B^*))$ is
dominating in $G$. Thus by the result of Step 1, the invariant subspace lattice
Lat$(T)$ for $T$  is rich, and so that $T$ has a infinitely many invariant
subspaces.
The proof is complete.\msk

\textbf{Acknowledgment.} The authors are deeply grateful to the referee for many valuable suggestions.

\end{document}